\documentclass[12pt]{article}
\usepackage{pifont}

\usepackage{amsfonts}
\usepackage{latexsym}
\usepackage{amsmath}
\usepackage{amssymb}
\usepackage{color}

\newtheorem{theorem}{Theorem}[section]

\newtheorem{proposition}[theorem]{Proposition}

\newtheorem{corollary}[theorem]{Corollary}

\newtheorem{example}[theorem]{Example}

\newenvironment{proof}[1][Proof]{\noindent \textbf{#1.} }{\ \ \  $\Box$}

\newtheorem{lemma}[theorem]{Lemma}
\newtheorem{definition}[theorem]{Definition}
\newtheorem{remark}[theorem]{Remark}

\title{Ergodicity of Invariant Capacity}

\date{}

\author{Chunrong Feng$^1$, Panyu Wu$^2$, Huaizhong Zhao$^1$\\
{\small $^1$Department of Mathematical Sciences, Loughborough University}\\{\small $^2$Zhongtai Securities Institute for Financial Studies, Shandong University}\\{\small C.Feng@lboro.ac.uk; wupanyu@sdu.edu.cn; H.Zhao@lboro.ac.uk}}

\begin{document}

\maketitle

\begin{abstract}
In this paper, we investigate capacity preserving transformations and their ergodicity.
We show that for any measurable transformation $\theta$ there always exists a $\theta$-invariant capacity. We investigate some limit properties under capacity spaces and then give the concept of ergodicity for a capacity preserving transformation.
Based on this definition, we give several characterizations of
ergodicity. In particular, we obtain a type of Birkhoff's ergodic theorem  and prove that the ergodicity of $\theta$ with respect to an upper probability is equivalent to the strong law of large numbers.\\
\par  $\textit{Keywords:}$  Capacity,  ergodicity, invariant set, invariant capacity, strong law of large numbers, Choquet integral.
\end{abstract}



\section{Introduction}\label{sec:intro}
In this paper, we investigate capacity preserving transformations and their ergodicity.
Capacity (or nonadditive probability) arises in modelling heterogeneous environments, for example, a financial market where biased beliefs of future price movements drive the decision of stock-market participants and create ambiguous volatility. It was pointed out that the additive probability theory is not adequate in either economics (see \cite{Gilboa2} and \cite{Schmeidler}) or statistics (see \cite{Walley1991}).
Dynamical systems on a capacity space concern transformations from the capacity space to itself. It is vitally important to study the dynamics of such transformations, of which little is known.
When additivity ceases to be valid, many classical results turn out to be invalid and the situations become more complicated.

The classical ergodic theory deals with a probability preserving map $\theta$ from $\Omega$ to $\Omega$ on a probability space $(\Omega,\mathcal{F},P)$. Let $\mathcal{G}$ denote all the invariant sets with respect to $\theta$. Then $\theta$ is called ergodic if any invariant set $B\in\mathcal{G}$ has either $P(B)=0$ or $P(B)=1$ which is equivalent to for any $B\in\mathcal{G}$, $P(B)=0$ or $P(B^c)=0$. This means that the dynamical system cannot be decomposed into different dynamical systems  (see \cite{Kifer} or \cite{Walters}). However, in the capacity theory, the equivalence is no longer true. To see this, let $(\mu,\bar{\mu})$ denote a pair of conjugate capacities on $\mathcal{F}$. Then
for any $ B\in\mathcal{G}$, $\mu(B) = 0$ or $1$ is equivalent to that for any  $ B\in\mathcal{G}$, $\mu(B) =0$ or $ \bar{\mu}(B^c)=0$. But it is not equivalent to either that for any  $ B\in\mathcal{G}$, $\mu(B) =0$ or $ \mu(B^c)=0$, or that for any  $ B\in\mathcal{G}$, $\bar{\mu}(B) =0$ or $ \bar{\mu}(B^c)=0$.
So how to define an ergodic transformation is an issue worthy of discussions.
Cerreia-Vioglio, Maccheroni and Marinacci called a capacity $\mu$ ergodic if $\mu (\mathcal{G}) = \{0, 1\}$
and then established an ergodic theorem for lower probabilities in \cite{cerreia}.

We do not think $\mu(\mathcal{G})=\{0,1\}$ is adequate for ergodicity in a capacity space as it is still possible that $\mu(B^c)=1$ when $\mu(B)=1$. Because if $\theta$ is ``ergodic'' in the sense of \cite{cerreia}, the space $\Omega$  may still be divided into two sets $B$ and $B^c$, each having full capacity but $\theta(B)=B$ and $\theta(B^c)=B^c$. In other words, $\theta$ is reducible while it is ``ergodic'' in the sense of \cite{cerreia} (see our Example \ref{ex3}).
The irreducibility condition is important as it is the essence of ergodicity.
There are also other papers attempting to investigate the ergodicity in capacity spaces or sublinear expectation spaces from different angles, one can see \cite{Cooman} and \cite{Huliwangzheng} and the references therein. None of these papers dealt with dynamical property of processes especially the non-decomposable property.

Inspired by the idea of \cite{fengzhao} for the ergodicity on sublinear expectation spaces, we add another condition to define the ergodicity in  the capacity space which is for any $\theta$-invariant set $B$, either $\mu(B)=0$ or $\mu(B^c)=0$.
Under this case, if $\theta$ is ergodic with respect to  capacity $\mu$ then the space $\Omega$ cannot be decomposed into two sets $B$ and $B^c$, each having positive capacity but $\theta(B)=B$ and $\theta(B^c)=B^c$. In other words, $\theta$ is irreducible.
We further obtain three equivalent characterizations of our ergodicity with respect to an upper probability: recurrence (Theorem \ref{th4}), the shift invariant random variable being a constant quasi-surely (Theorem \ref{th5}) and the time average of evolution of a random variable converging to a constant quasi-surely (Theorem \ref{th7}).

Though the ergodicity definition in \cite{cerreia} is not adequate, Cerreia-Vioglio, Maccheroni and Marinacci obtained an important result that for bounded random variable $\xi$, if the lower probability is $\theta$-invariant, then
$$v\left(\left\{\omega: \ \lim\limits_{n\to\infty}{1\over n}\sum_{k=0}^{n-1}\xi(\theta^{k}\omega) \hbox{ exists}\right\}\right)=1.$$
This is useful in our proof of the result that $\theta$ is ergodic with respect to an upper probability in the sense defined in this paper if and only if $\lim\limits_{n\to\infty}{1\over n}\sum_{k=0}^{n-1}\xi(\theta^{k}\omega)$ exists quasi-surely and is a constant quasi-surely.
Moreover, if the upper probability is continuous and concave, we show that this constant is bounded by the Choquet integrals with respect to the upper probability and the conjugate lower probability (Theorem \ref{th8}). This is due to the property that $\theta$ preserving capacity can infer $\theta$ preserving the corresponding Choquet integral (Proposition \ref{pr2}).

In the classical additive probability case, the result that any shift invariant random variable must be constant, as stated in Theorem \ref{th5} for capacities,
provides an important characteristics of ergodicity in terms of
the spectral structure of the corresponding transformation operator
on the space of measurable functions. For a Markov semigroup with an invariant measure, this
suggests that $0$ is a simple eigenvalue of the infinitesimal generator of the Markov semigroup
if and only if the invariant measure is ergodic. A well-known case in literature
is that of mixing stationary processes. In this case
the Koopman-von Neumann theorem implies that the generator has only one eigenvalue $0$ on the imaginary axis which is simple. Recently
it was observed in \cite{fengzhao2} that random periodic process is another ergodic regime of which the
spectral structure of the generator is distinct from that of the mixing regime.
In the ergodic random periodic regime
the infinitesimal generator
has infinite number of simple eigenvalues including $0$ equally placed on the imaginary axis.


The difference of the ergodicity of capacity developed in this paper and the one suggested by \cite{cerreia} can be demonstrated by looking at the well known Ellsberg Paradox (see \cite{ellsberg}). The Ellsberg Paradox says that there are two urns with red balls and black balls. In urn-I, both composition and proportion are known, say $50$ red balls and $50$ black balls. In urn-II, the proportion of red and black balls is not known. In the case of urn-I, it is logical to consider as equiprobable both the drawing of a red ball and that of a black ball. But considering of symmetry suggests that even for the urn-II,  the events of drawing a red ball and a black ball are still equal. However it is natural to expect that the confidence which the agent has in the two assignments is different. To reflect their lack of confidence and risk aversion, they may put
$$P(R_{II})=P(B_{II})=0.3 \hbox{ but still }P(R_{II}\cup B_{II})=1, $$
as their own probability judgement, which $R_{II}\ (B_{II})$ is the event of drawing a red (black) ball from the urn-II. This typical case of subjective probability is an example of capacity. Now one can place countable number of identical urn-IIs and draw a ball from each of these urns. Define $X_n=0$ when drawing a black ball from the $n$-th urn-II and  $X_n=1$ when drawing a red ball from the $n$-th urn-II. Certainly, despite the nonadditive subjective probability, an agent has, in reality, as $n\to\infty$, the empirical average
$\frac{1}{n}\sum_{i=1}^n X_i$ converges to the true proportion of red balls in urn-II, which is the objective probability. The agent does not know the objective probability in the first place, nevertheless the convergence to the objective probability does not depend on the agent's subjective probability. However the theory developed in \cite{cerreia} does not imply $\frac{1}{n}\sum_{i=1}^n X_i$ converges to the objective probability. In stead, their theory can only imply the  empirical mean converges but the limit is a random variable and cannot be pinned down to  a constant. In their definition of ergodicity, there might exist an invariant set that has zero lower probability but is still quite ``large''. This reducibility stops one concluding that $\frac{1}{n}\sum_{i=1}^n X_i$ converges to a constant. The theory we develop in this paper say that the necessary and sufficient conditions for $\frac{1}{n}\sum_{i=1}^n X_i$ converges to a constant is the ergodicity of upper probability in the sense introduced in this paper. When this result is applied to the Ellsberg box model,  it reveals the reality that the limit should be exactly the true proportion of red balls. This is made possible as our definition avoids the case of having a ``large'' invariant set with zero capacitliy.

So far, no other limit theory about capacity provided a condition for
the limit of $\frac{1}{n}\sum_{i=1}^n X_i$ being a constant. The ergodic theory of this paper can provide new insight to the study of capacity. Needless to say it  can apply to many other problems apart from the Ellsberg Paradox, of which the limit theory is obvious.

The Kolmogorov $0$-$1$ Law plays an important role in the limit theory under
the classical probability framework, which implies that the tail event happens with probability $0$ or $1$.
In Section 2, we investigate the Kolmogorov $0$-$1$ Law in the capacity space setting.
We give an example to show that a sequence of independent random variables with respect to capacity $\mu$ cannot deduce
$\sigma(Y_k, k \le n)$ and $\sigma(Y_k, k \ge n + 1)$ being independent with respect to  $\mu$, for any
$n \in\mathbb{N}$. Then we investigate some properties of Choquet integral which will be used in this paper.
In Section 3, we study the properties of $\theta$-invariant capacities. We show that for any measurable transformation $\theta$ there always exists a   $\theta$-invariant capacity. That is a surprising result which is not true in the classical probability case.
In Section 4, we firstly investigate the limit properties under the scenario that $\mu(\mathcal{G})=\{0,1\}$ and then give the definition of a  transformation to be ergodic.
Based on our definition, we give several characterisations of
ergodicity and a type of Birkhoff's ergodic theorem. In Section 5, we give a strong law of large numbers for stationary and ergodic sequences in an upper probability space.

\section{Basic concepts and independence on capacity space}
Let $(\Omega,\mathcal{F})$ be a measurable space. Recall a set function
$\mu:\mathcal{F}\rightarrow[0, 1]$ is

(i) a capacity/nonadditive probability if $\mu(\emptyset)=0,$ $\mu(\Omega)=1$, and $\mu(A)\leq\mu(B)$ for all $A, B \in\mathcal{F}$ such that $A\subseteq B$;

(ii) concave/submodular/2-alternating if $\mu(A\cup B)+\mu(A\cap B)\leq\mu(A)+\mu(B)$ for all $A, B \in\mathcal{F}$;

(iii) convex/supermodular/2-monotonic if $\mu(A\cup B)+\mu(A\cap B)\geq\mu(A)+\mu(B)$ for all $A, B \in\mathcal{F}$;

(iv) subadditive if $\mu(A\cup B)\leq\mu(A)+\mu(B)$ for all $A, B \in\mathcal{F}$ with $A\cap B =\emptyset$;

(v) superdditive if $\mu(A\cup B)\geq\mu(A)+\mu(B)$ for all $A, B \in\mathcal{F}$ with $A\cap B =\emptyset$;

(vi) continuous from below/inner continuous if $\lim\limits_{n\to\infty}\mu(A_n)=\mu(A)$ for $A_n\uparrow A$;

(vii) continuous from above/outer continuous if $\lim\limits_{n\to\infty}\mu(A_n)=\mu(A)$ for $A_n\downarrow A$;

(viii) continuous if it is both continuous from below and above;

(ix) continuous at $\emptyset$ if $\lim\limits_{n\to\infty}\mu(A_n)=0$ for $A_n\downarrow\emptyset$;

(x) continuous at $\Omega$ if $\lim\limits_{n\to\infty}\mu(A_n)=1$ for $A_n\uparrow\Omega$.

It is obvious that $\mu$ being concave/convex implies $\mu$ being subadditive/superadditive. The converse is not correct. It is also easy to check that for a subadditive capacity $\mu$, $\mu$ being continuous from above is equivalent to $\mu$ being continuous at $\emptyset$.

For a capacity $\mu$ on $\mathcal{F}$, we call $(\Omega,\mathcal{F},\mu)$ a capacity space.
We can define the
conjugate capacity $\bar{\mu}$ on $\mathcal{F}$ by $$\bar{\mu}(A)=1-\mu({A^c}),\qquad \hbox{for any } A\in\mathcal{F},$$
where $A^c$ is the complementary set of $A$.
Notice that if $\mu$ is additive, then $\bar{\mu} = \mu$. Capacity $\mu$ is continuous at $\Omega$ if and only if $\bar{\mu}$ is continuous at $\emptyset$.
Moreover, capacity $\mu$ is convex if and only if $\bar{\mu}$ is concave. However, the superadditivity and the subadditivity do not have such a conjugate relation.

Let $\Delta(\Omega,\mathcal{F})$ denote the set of all finitely additive probabilities on $\mathcal{F}$ and $\Delta^{\sigma}(\Omega,\mathcal{F})$ denote the set of all  probabilities ($\sigma$-additive) on $\mathcal{F}$.
The widely studied conjugate capacities which satisfy the subadditivity and superadditivity respectively are upper and lower probabilities. A pair of capacities $(V,v)$ is called the upper and lower probabilities on $(\Omega,\mathcal{F})$ (generated by $\mathcal{P}$) if
$$V(A)=\sup _{P\in\mathcal{P}}P(A) \hbox{ and } v(A)=\inf _{P\in\mathcal{P}}P(A), \hbox{ for any } A\in\mathcal{F},$$
where $\mathcal{P}$ is a nonempty set of $\Delta(\Omega,\mathcal{F})$. The core of lower probability $v$ is defined by
$$\hbox{core}(v)=\{P\in\Delta(\Omega,\mathcal{F}):\ P\ge v\}.$$
In the sequel, we use $\mathbb{N}$ to denote the set of all the positive integers and $\mathbb{N}_0=\mathbb{N}\cup\{0\}$.

\begin{proposition}\label{pr5}
Let $v$ be a lower probability on $(\Omega,\mathcal{F})$ generated by $\mathcal{P}$, where $\mathcal{P}$ is a nonempty set of $\Delta(\Omega,\mathcal{F})$. Then \\
(i) $\hbox{core}(v)\neq\emptyset$ and
$\hbox{core}(v)$ is weak$^{*}$ compact;\\
(ii) lower probability $v$ is exact, that is to say for any $A\in\mathcal{F}$, $v(A)=\min _{P\in \hbox{core}(v)}P(A).$
\end{proposition}
\begin{proof}
(i) It is obvious that $\hbox{core}(v)\neq\emptyset$ since $\mathcal{P}\subseteq \hbox{core}(v)$ and the weak$^{*}$  compactness of
$\hbox{core}(v)$ is directly deduced from Proposition 4.2 in \cite{Marinacci3}.

(ii) For any $A\in\mathcal{F}$, on one hand,
$ v(A)=\inf _{P\in\mathcal{P}}P(A)\ge \inf _{P\in \hbox{core}(v)}P(A)$
since $\mathcal{P}\subseteq \hbox{core}(v)$. On the other hand, by the definition of $\hbox{core}(v)$, it is clear that  $ v(A)\le \inf _{P\in \hbox{core}(v)}P(A)$.
Therefore $$ v(A)=\inf _{P\in \hbox{core}(v)}P(A), \hbox{ for any } A\in\mathcal{F}. $$
Thus, for any given $A\in\mathcal{F}$, $n\in \mathbb{N}$, there exists $P_n\in\hbox{core}(v)$ such that $P_n(A)\le v(A)+{1\over n}$.
Due to the weak$^{*}$  compactness of
$\hbox{core}(v)$, there exists a subsequence $\{P_{n_k}\}_{k=1}^{\infty}\subseteq\{P_n\}_{n=1}^{\infty}$ such that $P_{n_k}$ weak$^{*}$ converges to a finitely additive probability $\bar{P}$ and $\bar{P}\in \hbox{core}(v)$. Hence,
$$\bar{P}(A)=\lim_{k\to\infty}P_{n_k}(A)\le v(A),$$
as well as $\bar{P}(A)\ge v(A)$. Therefore $v(A)=\min _{P\in \hbox{core}(v)}P(A).$
\end{proof}

%
%

\begin{proposition}\label{pr6}
Let $(V,v)$ be a pair of upper and lower probabilities on $(\Omega,\mathcal{F})$ generated by $\mathcal{P}$, where $\mathcal{P}$ is a nonempty set of $\Delta(\Omega,\mathcal{F})$. Then the following conditions are equivalent \\
(i) upper probability $V$ is continuous at $\emptyset$;\\
(ii)  upper probability $V$ is continuous;\\
(iii) lower probability $v$ is continuous;\\
(iv) lower probability $v$ is continuous at $\Omega$.\\
And any statement of (i)-(iv) implies\\
(v) $\hbox{core}(v)$ is a subset of $\Delta^{\sigma}(\Omega,\mathcal{F})$.
\end{proposition}
\begin{proof}
It is easy to check that (i)$\Leftarrow$(ii)$\Leftrightarrow$(iii)$\Rightarrow$(iv)$\Leftrightarrow$(i). So we only need to prove (i)$\Rightarrow$(ii) and (i)$\Rightarrow$(v).

Due to Proposition \ref{pr5}, we have $\hbox{core}(v)\neq\emptyset$.
Fix any $P\in \hbox{core}(v)$, it follows from (i) that $P$ is $\sigma$-additive, since
$$\limsup_{n\to\infty}P(A_n)\le \lim_{n\to\infty}V(A_n)=0,\quad\hbox{for all } A_n\downarrow\emptyset.$$
In other words, statement (v) holds. Then it is easy to check that $V$ is continuous from below (or see Lemma 2.1 in \cite{chenwuli}).

For any $A_n\downarrow A$, $A_n,A\in\mathcal{F}$, 
we have $V(A)\le V(A_n)$ for all $n\ge1$.
By $A_n\setminus A \downarrow\emptyset$, subadditivity of $V$ and (i), we have
$$\lim_{n\to\infty}V(A_n)\le \lim_{n\to\infty}V(A_n\setminus A)+V(A)=V(A).$$
Therefore, $V(A)=\lim_{n\to\infty}V(A_n)$, that is, $V$ is continuous from above. Hence, statement (ii) holds.
\end{proof}

\begin{remark}
If $\Omega$ is a  Polish space and $\mathcal{F}$ is the Borel $\sigma$-algebra, then the statements (i)-(iv) are equivalent to $\hbox{core}(v)$ being a weakly compact subset of $\Delta^{\sigma}(\Omega,\mathcal{F})$. This result can be obtained from the above Proposition \ref{pr6}  and Theorem 4.2 in \cite{Marinacci3}.
\end{remark}

Before we establish the Kolmogorov $0$-$1$ Law in a capacity space, we give the following definitions which are natural extensions of the corresponding concepts in the classical probability theory.
\begin{definition}\label{de4}
Let $(\Omega,\mathcal{F},\mu)$ be a  capacity space and $J$ be an index set.

Events $A$ and $B$ are called independent with respect to $\mu$ if $\mu(A \cap B) = \mu(A) \mu(B )$.

Events set $\{A_t,\ t\in J\}$ are called (mutually) independent with respect to $\mu$ if for any finite subset $I\subseteq J$
$$\mu\left(\bigcap\limits_{t\in I}A_t\right) = \prod\limits_{t\in I}\mu(A_t).$$

Let $\{\mathcal{D}_t,\ t\in J\}$ be a group of subclasses of $\mathcal{F}$. If for any finite subset $I\subseteq J$
$$\mu\left(\bigcap\limits_{t\in I}A_t\right) = \prod\limits_{t\in I}\mu(A_t),\quad \hbox{for any } A_t\in\mathcal{D}_t,\  t\in I,$$
then $\{\mathcal{D}_t,\ t\in J\}$ are called (mutually) independent subclasses with respect to  $\mu$.
\end{definition}

\begin{lemma}\label{le4}
If $\{\mathcal{D}_t,\ t\in J\}$ are independent subclasses with respect to  a continuous capacity $\mu$ and $\mathcal{D}_t$ is an algebra for every $t\in J$, then $\{\sigma(\mathcal{D}_t),\ t\in J\}$ are independent $\sigma$-algebra with respect to  the capacity $\mu$.
\end{lemma}
\begin{proof}
This lemma can be deduced from the monotone class theorem. We omit the details.
\end{proof}

\begin{definition}\label{de5}
Random variables $\{Y_t,t\in J\}$ on $(\Omega,\mathcal{F},\mu)$ are said to be independent with respect to  $\mu$ if the $\sigma$-algebra
$\{\sigma(Y_t),t\in J\}$ are independent with respect to  $\mu$.
\end{definition}

For a sequence of random variables $\{Y_n\}_{n\in\mathbb{N}}$, Definition \ref{de5} is equivalent to the definition of independent random variables sequence with respect to  $\mu$ given by \cite{Marinacci2} (Definition 4). But it is worth noting that $\{Y_n\}_{n\in\mathbb{N}}$ being independent with respect to  $\mu$ cannot deduce $\sigma(Y_k,\ k\le n)$ and  $\sigma(Y_k,\ k\ge n+1)$ being independent with respect to  $\mu$, for any $n\in\mathbb{N}$. We give the following example to illustrate this.

\begin{example}\label{ex2}
Let $\Omega=\{\omega_1,\omega_2,\omega_3,\omega_4,\omega_5\}$, $\mathcal{F}$ be all subsets of $\Omega$ and $\mu$ be a capacity on  $\mathcal{F}$, with $\mu(A)=0$ if $|A|\le 3 $, and $\mu(A)=1$ if $|A|\ge4$, where $|A|$ denotes the  number of elements in $A$. It is easy to see that $\mu$ is superadditive but is not convex. However, the conjugate capacity of $\mu$ is neither subadditive nor superadditive.
Let
\begin{equation*}
Y_1(\omega)=
\begin{cases}
0& \omega=\omega_1,\omega_2\\
1& \omega=\omega_3,\omega_4,\omega_5
\end{cases}
,\quad
Y_2(\omega)=
\begin{cases}
0& \omega=\omega_1,\omega_4,\omega_5\\
1& \omega=\omega_2,\omega_3
\end{cases}
,
\end{equation*}
\begin{equation*}
Y_3(\omega)=
\begin{cases}
0& \omega=\omega_3,\omega_4\\
1& \omega=\omega_1,\omega_2,\omega_5
\end{cases}
,\quad
Y_4(\omega)=
\begin{cases}
0& \omega=\omega_1,\omega_4\\
1& \omega=\omega_2,\omega_3,\omega_5
\end{cases}
.
\end{equation*}
It is easy to check $Y_1,Y_2,Y_3,Y_4$ are independent and identically distributed with respect to  $\mu$. However, $\sigma(Y_1,Y_2)$ and $\sigma(Y_3,Y_4)$ are not independent  with respect to  $\mu$ since
\begin{eqnarray*}&&\mu((\{Y_1=0\}\cup\{Y_2=0\})\cap(\{Y_3=0\}\cup\{Y_4=1\}))=\mu(\{\omega_2,\omega_4,\omega_5\})=0\\
&&\neq 1 =\mu(\{Y_1=0\}\cup\{Y_2=0\})\cdot\mu(\{Y_3=0\}\cup\{Y_4=1\})).
\end{eqnarray*}

Moreover, notice that
$$\mu(\omega:\ (Y_1,Y_2)(\omega)\in\left\{(0,0),(0,1),(1,1)\right\})=\mu(\omega_1,\omega_2,\omega_3)=0$$
but
$$\mu(\omega:\ (Y_3,Y_4)(\omega)\in\left\{(0,0),(0,1),(1,1)\right\})=\mu(\omega_2,\omega_3,\omega_4,\omega_5)=1.$$
So $(Y_1,Y_2)$ and $(Y_3,Y_4)$ are not identically distributed with respect to  $\mu$.
\end{example}

Next we give the Kolmogorov $0$-$1$ Law in capacity spaces that will be used in Section 5.

\begin{theorem}\label{th3}
Let $\{Y_n\}_{n\in\mathbb{N}}$ be random variables such that for any $n\in\mathbb{N}$, $\sigma(Y_k,\ k\le n)$ and  $\sigma(Y_k,\ k\ge n+1)$ are independent with respect to  a continuous capacity $\mu$. By $\mathcal{T}$ we denote the tail $\sigma$-algebra of $\{Y_n\}_{n\in\mathbb{N}}$, that is $\mathcal{T}=\bigcap\limits_{n=1}^{\infty}\sigma(Y_k,\ k\ge n)$.
Then for any $A\in \mathcal{T}$, the following two statements are true:\\
(i) $\mu(A)=0$ or $1$;\\
(ii) $\mu(A)=0$ or $\mu(A^c)=0$.
\end{theorem}
\begin{proof}
It is obvious that $\mathcal{T}$ is independent of $\mathcal{A}=\bigcup
_{n\ge1}\sigma(Y_k,\ k\le n)$ with respect to   $\mu$
since $\mathcal{T}\subseteq \sigma(Y_k,\ k\ge n+1)$, for all $n\in \mathbb{N}$.
Notice that $\mathcal{A}$ is an algebra, then by Lemma \ref{le4}, $\mathcal{T}$ is independent of $\sigma(\mathcal{A})$. On the other hand, $\mathcal{T}\subseteq\sigma(\mathcal{A})$.
Hence, $\mathcal{T}$ is independent of itself  with respect to   $\mu$. So for any $A\in \mathcal{T}$, $\mu(A)=\mu(A\cap A)=\mu(A)\cdot \mu(A)$, which implies $\mu(A)=0$ or $1$. Meantime, $\mu(\emptyset)=\mu(A\cap A^c)=\mu(A)\cdot \mu(A^c)$, which deduces $\mu(A)=0$ or $\mu(A^c)=0$.
\end{proof}

\begin{remark}\label{re1}
Although the proof of the Kolmogorov $0$-$1$ Law in capacity spaces here is no big difference with the  proof of the classical Kolmogorov $0$-$1$ Law, we still present here since the two statements (i) and (ii) in Theorem \ref{th3} are  equivalent when $\mu$ is additive, but not equivalent when $\mu$ is nonadditive.
\end{remark}

At the end of this section, we recall the Choquet integral/expectation of a random variable, introduced by Choquet in \cite{Choquet}.
For any $\mathcal{F}$-measurable real valued random variable $\xi$, the Choquet integral/expectation of $\xi$ with respect to  $\mu$ is defined by
$$\int_{\Omega} \xi(\omega) d\mu=\int_{0}^{\infty}\mu(\{\omega:\ \xi(\omega)\geq t\})\mathrm{d}t+\int_{-\infty}^0[\mu(\{\omega:\ \xi(\omega)\geq t\})-1]\mathrm{d}t.$$ In this paper, we always consider the random variables taking real values. The asymmetry is one of the most important properties of Choquet integral (see Proposition 5.1 in \cite{Denneberg}), which means that
$$\int_{\Omega}-\xi(\omega)d\mu=-\int_{\Omega}\xi(\omega)d\bar{\mu}.$$
It is well known that the Choquet integral with respect to $\mu$ will become a sublinear expectation (see \cite{p2010}) when $\mu$ is concave. We put this property here as Lemma \ref{le5} since it will the used in the sequel.

\begin{lemma}\label{le5}
Let $\mu$ be a concave capacity on $\mathcal{F}$. Then for random variables $\xi_i$ with $\int_{\Omega}\xi_i(\omega)d\mu>-\infty$, $i=1,2$, we have
$$\int_{\Omega}(\xi_1+\xi_2)(\omega)d\mu\le\int_{\Omega}\xi_1(\omega)d\mu+\int_{\Omega}\xi_2(\omega)d\mu.$$
If $\mu$ is continuous from below the assumption on $\xi_1$ and $\xi_2$ can be dropped.
\end{lemma}
\begin{proof}
This lemma is a direct corollary of Theorem 6.3 in \cite{Denneberg}.
\end{proof}

If further $\mu$ is continuous at $\emptyset$, then the Choquet integral with respect to $\mu$ will be an upper expectation as the following Proposition \ref{pr4} shows.

\begin{proposition}\label{pr4}
Let $\mu$ be a concave capacity on $\mathcal{F}$ and be continuous at $\emptyset$. Then\\
(i) there exists a set of $\sigma$-additive probabilities $\mathcal{P}$ such that  $$\mu(A)=\sup_{P\in\mathcal{P}}P(A),\quad \hbox{for all }A\in\mathcal{F},$$
and moreover,
$$\int_{\Omega}\xi(\omega)d\mu=\sup_{P\in\mathcal{P}}\int_{\Omega}\xi(\omega)dP,\quad \hbox{for all } \xi \hbox{ such that } \int_{\Omega}|\xi(\omega)|d\mu<\infty;$$
(ii) capacity $\mu$ is continuous.
\end{proposition}
\begin{proof}
Since $\mu$ is concave, by Proposition 10.3 in \cite{Denneberg}, there exists a nonempty set
$$\mathcal{M}=\{\alpha:\ \alpha \hbox{ is additive capacity and } \alpha(A)\le\mu(A),\ \hbox{for all } A\in\mathcal{F}   \}$$
such that
$$\int_{\Omega}\xi(\omega)d\mu=\sup_{\alpha\in\mathcal{M}}\int_{\Omega}\xi(\omega)d\alpha,\quad \hbox{for all } \xi \hbox{ such that } \int_{\Omega}|\xi(\omega)|d\mu<\infty.$$
So $\mu$ is an upper probability. Then we can prove (i) and (ii) by Proposition \ref{pr6}.
%
\end{proof}

The following lemma can be seen as the dominated convergence theorem in a capacity space with respect to the Choquet integral and will be used in the proof of Theorem \ref{th8}.

\begin{lemma}\label{le2}
Let $\mu$ be a subadditive capacity on $\mathcal{F}$ and is continuous from above, $\{X_n\}_{n\in\mathbb{N}}$, $Y$ and $Z$ be random variables with $Y\le X_n\le Z$
and $\int_\Omega Y d\mu$, $\int_\Omega Z d\mu$ finite. If $\bar{\mu}\left(\{\omega:\ \lim\limits_{n\to\infty}X_n(\omega)= X(\omega)\}\right)=1$, then
\begin{equation}\label{eq16}
\lim_{n\to\infty}\int_\Omega X_n d\mu=\int_\Omega X d\mu.
\end{equation}
\end{lemma}
\begin{proof}
For any $\epsilon >0$, we get
$$\bar{\mu}\left(\bigcup_{n=1}^{\infty}\bigcap_{k=n}^{\infty}\left\{\omega:\ |X_k(\omega)- X(\omega)|<\epsilon\right\}\right)=1$$
from
\begin{eqnarray*}
\bar{\mu}\left( \bigcup_{n=1}^{\infty}\bigcap_{k=n}^{\infty}\left\{\omega:\  |X_k(\omega)- X(\omega)|<\epsilon\right\}\right)
&\ge& \bar{\mu}\left( \bigcap_{\epsilon >0}\bigcup_{n=1}^{\infty}\bigcap_{k=n}^{\infty} \{\omega:\ |X_k(\omega)- X(\omega)|<\epsilon\}\right)\\
&=&\bar{\mu}\left(\{\omega:\ \lim_{n\to\infty}X_n(\omega)= X(\omega)\}\right)=1.
\end{eqnarray*}
Therefore, by the continuity from above of $\mu$, we have
\begin{eqnarray*}
0&\le&\limsup\limits_{n\to\infty} \mu\left(\{\omega:\ |X_n(\omega)- X(\omega)|\ge\epsilon\}\right)\\
&\le&\limsup\limits_{n\to\infty} \mu\left(\bigcup_{k=n}^{\infty}\left\{\omega:\ |X_k(\omega)- X(\omega)|\ge\epsilon\right\}\right)\\
&=&\mu\left(\bigcap_{n=1}^{\infty}\bigcup_{k=n}^{\infty}\left\{\omega:\ |X_k(\omega)- X(\omega)|\ge\epsilon\right\}\right)\\
&=&1-\bar{\mu}\left(\bigcup_{n=1}^{\infty}\bigcap_{k=n}^{\infty}\{|X_k(\omega)- X(\omega)|<\epsilon\}\right)=0,
\end{eqnarray*}
which means $X_n$ converges to $X$ $\mu$-stochastically (p97 in \cite{Denneberg}). Since $\mu$ is subadditive, then by Proposition 8.5  and Theorem 8.9 in \cite{Denneberg}, we can get (\ref{eq16}).
\end{proof}

\section{Invariant capacity and its existence}
We consider a $\mathcal{F}/\mathcal{F}$-measurable transformation $\theta:\ \Omega\to\Omega$.
A set $A$ is called invariant set with respect to  $\theta$ if $\theta^{-1}A=A$. It is easy to check that $A^c$ is an invariant set if and only if $A$ is an invariant set. Let $\mathcal{G}$ denote the set of all invariant sets with respect to  $\theta$, it is easy to check that $\mathcal{G}$ is a sub-$\sigma$-algebra of $\mathcal{F}$ (Exercise 7.1.1 in \cite{Durrett}).
Corresponding to the $\theta$-invariant probability, Cerreia-Vioglio, Maccheroni and Marinacci introduced the definition of $\theta$-invariant capacity (Definition 1 in \cite{cerreia}) as follows.

\begin{definition}\label{de1}
A capacity $\mu$ is $\theta$-invariant if for each $A\in\mathcal{F}$, $\mu(A)=\mu(\theta^{-1}A).$
We also say that $\theta$ preserves $\mu$ if $\mu$ is $\theta$-invariant.
\end{definition}

Firstly, we will show for any $\mathcal{F}/\mathcal{F}$-measurable transformation $\theta$, there exists a $\theta$-invariant capacity.
This property actually is not enjoyed  in the additive probability case, which is the main challenge of the classical theory. Let $\mathcal{N}$ denote the set of all $\theta$-invariant capacities, $\overline{\mathcal{N}}$ denote the set of all subadditive $\theta$-invariant capacities and $\underline{\mathcal{N}}$ denote the set of all superadditive $\theta$-invariant capacities.

\begin{proposition}\label{pr1}
For any $\mathcal{F}/\mathcal{F}$-measurable transformation $\theta:\ \Omega\to\Omega$, $\overline{\mathcal{N}}\neq\emptyset$ and $\underline{\mathcal{N}}\neq\emptyset$. In particular, $\mathcal{N}\neq\emptyset$.
\end{proposition}
\begin{proof}
Fix an arbitrary $\omega\in\Omega$, for any $A\in\mathcal{F}$, set
$$\mu(A)=\limsup\limits_{n\to\infty}{1\over n}\sum\limits_{k=0}^{n-1}I_{A}(\theta^{k}\omega)=\limsup\limits_{n\to\infty}{1\over n}\sum\limits_{k=0}^{n-1}I_{\theta^{-k}A}(\omega),$$
which can be seen as the super visit frequency of $\theta^{k}\omega$ to the set $A$. It is easy to check that
$\mu$ is a subadditive capacity on $\mathcal{F}$. Next we prove $\mu$ is $\theta$-invariant. For this,
\begin{eqnarray*}
\left|\mu(\theta^{-1}A)-\mu(A)\right|&=&\left|\limsup\limits_{n\to\infty}{1\over n}\sum\limits_{k=0}^{n-1}I_{\theta^{-k-1}A}(\omega)-\limsup\limits_{n\to\infty}{1\over n}\sum\limits_{k=0}^{n-1}I_{\theta^{-k}A}(\omega)\right|\\
&\le&\limsup\limits_{n\to\infty}\left|{1\over n}\sum_{k=0}^{n-1}[I_{\theta^{-k-1}A}(\omega)-I_{\theta^{-k}A}(\omega)]\right|\\
&=&\limsup\limits_{n\to\infty}{1\over n}\left|I_{\theta^{-n}A}(\omega)-I_{A}(\omega)\right|\\
&\le&\limsup\limits_{n\to\infty}{2\over n}=0.
\end{eqnarray*}
Hence
$\mu(\theta^{-1}A)=\mu(A)$, that is $\mu\in\overline{\mathcal{N}}$.
Then the conjugate capacity $\bar{\mu}$ satisfies $\bar{\mu}(A)=\liminf\limits_{n\to\infty}{1\over n}\sum\limits_{k=0}^{n-1}I_{A}(\theta^{k}\omega)$, thus $\bar{\mu}$ is supperadditive and $\bar{\mu}\in\underline{\mathcal{N}}$. Finally the claim $\mathcal{N}\neq\emptyset$ follows from $\overline{\mathcal{N}}$ and $\underline{\mathcal{N}}$ being nonempty sets.
\end{proof}

\begin{proposition}\label{pr2}
(i) A capacity $\mu$ is $\theta$-invariant if and only if its conjugate capacity $\bar{\mu}$ is $\theta$-invariant.

(ii) The sets $\mathcal{N}$, $\overline{\mathcal{N}}$ and $\underline{\mathcal{N}}$ are convex.

(iii) If the capacity $\mu$ is $\theta$-invariant, then $\theta$ preserves the Choquet integral with respect to  $\mu$, that is
$$\int_{\Omega} \xi(\theta\omega) d\mu=\int_{\Omega} \xi(\omega) d\mu.$$
\end{proposition}

\begin{proof} (i) and (ii) are easy to check by the definition of $\theta$-invariant and conjugate capacity.

(iii) To check the invariance of the Choquet integral, we consider
\begin{eqnarray*}
&&\int_{\Omega} \xi(\theta\omega) d\mu
=\int_{0}^{\infty}\mu(\{\omega:\ \xi(\theta\omega)\geq t\})\mathrm{d}t+\int_{-\infty}^0[\mu(\{\omega:\ \xi(\theta\omega)\geq t\})-1]\mathrm{d}t\\
&=&\int_{0}^{\infty}\mu(\{\omega:\ \omega\in\theta^{-1}(\xi^{-1}[t,\infty))\})\mathrm{d}t+\int_{-\infty}^0[\mu(\{\omega:\ \omega\in\theta^{-1}(\xi^{-1}[t,\infty))\})-1]\mathrm{d}t\\
&=&\int_{0}^{\infty}\mu(\{\omega:\ \omega\in\xi^{-1}[t,\infty)\})\mathrm{d}t+\int_{-\infty}^0[\mu(\{\omega:\ \omega\in\xi^{-1}[t,\infty)\})-1]\mathrm{d}t\\
&=&\int_{\Omega} \xi(\omega) d\mu,
\end{eqnarray*}
where the penultimate equality was due to the $\theta$-invariance of $\mu$.
\end{proof}

The following lemma can be deduced from the Theorem 2 obtained by Cerreia-Vioglio, Maccheroni and Marinacci  in \cite{cerreia} which will be useful in our characterization about the ergodicity in an upper probability space.

\begin{lemma}\label{pr7}
Let $v$ be a continuous lower probability on $(\Omega,\mathcal{F})$. If $v$ is  $\theta$-invariant, then for any bounded $\mathcal{F}$-measurable random variable $\xi$,
$$v\left(\left\{\omega:\ \lim_{n\to\infty}{1\over n}\sum_{k=0}^{n-1}\xi(\theta^{k}(\omega)) \hbox{ exists }\right\}\right)=1.$$
\end{lemma}

\section{Ergodicity under capacity space}
Before studying the ergodicity under a capacity space, we firstly study the following properties of random variables which are measurable with respect to a sub-$\sigma$-algebra of $\mathcal{F}$ with capacity only $0$ or $1$.

\begin{theorem}\label{le1}
Let $\mu$ be a continuous capacity on $\mathcal{F}$ and $\mathcal{F}_0$ be any sub-$\sigma$-algebra of $\mathcal{F}$ with $\mu(\mathcal{F}_0)=\{0,1\}$. Then

(i) for any $\mathcal{F}_0$-measurable random variable $\xi$, we have
\begin{equation}\label{eq1}
\mu\left\{\omega:\  \xi(\omega)\ge \int_{\Omega}\xi d\mu\right\}=1,
\end{equation}
and
\begin{equation}\label{eq2}
\mu\left\{\omega:\ \xi(\omega)\le\int_{\Omega}\xi d\bar{\mu}\right\}=1;
\end{equation}

(ii) if $\mu(A\cap B)=1$ for any $A,\ B\in\mathcal{F}_0$ with $\mu(A)=1$ and $\mu(B)=1$, then for any $\mathcal{F}_0$-measurable random variables $\xi$, we have
\begin{equation}\label{eq3}
\mu\left\{\omega:\ \int_{\Omega}\xi d\mu \le\xi(\omega)\le\int_{\Omega}\xi d\bar{\mu}\right\}=1.
\end{equation}
\end{theorem}

\begin{proof} (i) Since $\mu$ is continuous from below, we have
\begin{eqnarray*}1=\mu\left(\left\{\omega:\ \xi(\omega)\in(-\infty,\infty)\right\}\right)&=&\mu\left(\bigcup\limits_{n=1}^{\infty}\left\{\omega:\ \xi(\omega)\in[-n,n]\right\}\right)\\
&=&\lim\limits_{n\to\infty}\mu\left(\{\omega:\ \xi(\omega)\in [-n,n]\}\right).
\end{eqnarray*}
Notice that $\mu(\mathcal{F}_0)=\{0,1\}$ and $\xi$ is $\mathcal{F}_0$-measurable, therefore there exists $n\in\mathbb{N}$, such that $\mu\left(\{\omega:\ \xi(\omega)\in [-n,n]\}\right)=1$. It turns out that the following set $I$ is not empty, where
$$I=\{t\in \mathbb{R}:\ \mu\{\omega:\ \xi(\omega)\ge t\}=1\}.$$
We define
$t^*=\sup I$.
Since $\mu$ is continuous from above, we get $\mu\left(\{\omega:\ \xi(\omega)\ge t^*\}\right)=1$, so $t^*\in I$.
Due to $\mu(\mathcal{F}_0)=\{0,1\}$, for any $t>t^*$, $\mu\left(\{\omega:\ \xi(\omega)\ge t\}\right)=0$, and for any $t\le t^*$, $\mu\left(\{\omega:\ \xi(\omega)\ge t\}\right)=1$.
The above observations lead to that, if $t^*\ge0$ then
\begin{eqnarray*}
\int_{\Omega}\xi d\mu=\int_{0}^{\infty}\mu\left(\{\omega: \ \xi(\omega)\ge t\}\right)\mathrm{d}t+\int_{-\infty}^0[\mu(\{\omega:\ \xi(\omega)\geq t\})-1]\mathrm{d}t
=\int_{0}^{t^{*}}1\mathrm{d}t=t^{*};
\end{eqnarray*}
and if $t^*<0$ then
\begin{eqnarray*}
\int_{\Omega}\xi d\mu=\int_{0}^{\infty}\mu\left(\{\omega: \ \xi(\omega)\ge t\}\right)\mathrm{d}t+\int_{-\infty}^0[\mu(\{\omega:\ \xi(\omega)\geq t\})-1]\mathrm{d}t=\int^{0}_{t^{*}}(-1)\mathrm{d}t=t^{*}.
\end{eqnarray*}
Therefore $t^{*}=\int_{\Omega}\xi d\mu$ and the equality (\ref{eq1}) holds.

Considering random variable $-\xi$ in (\ref{eq1}), we get $\mu\left(\{\omega:\ -\xi(\omega)\ge \int_{\Omega}-\xi d\mu \}\right)=1$.
Thus we obtain that the equality (\ref{eq2}) holds since $\int_{\Omega}-\xi d\mu=-\int_{\Omega}\xi d\bar{\mu}$.

(ii) Under the new assumption, equality (\ref{eq3}) can be deduced  directly from equalities (\ref{eq1}) and (\ref{eq2}).
\end{proof}

\begin{remark}
(i) If $\xi$ is bounded, then we can replace the continuity of $\mu$ by the continuity from below of $\mu$, the conclusions in Theorem \ref{le1} still hold.

(ii) From the proof of Theorem \ref{le1} we can see that if $\mu$ is continuous and $\mu(\mathcal{F}_0)=\{0,1\}$, then for any $\mathcal{F}_0$-measurable random variable $\xi$, $\xi$ is bounded $\mu$-a.e.

(iii) When $\xi$ is bounded and $\mu$ is lower probability, then Theorem \ref{le1} degenerates to Lemma 2 in \cite{cerreia}.
\end{remark}

If $\mu$ is a lower probability, it satisfies that $\mu(A\cap B)=1$ for any $\mu(A)=1$ and $\mu(B)=1$. But there are also other capacities rather than lower probabilities satisfying this condition. See the following example.

\begin{example}
Let $(\Omega,\mathcal{F},P)$ be a probability space and $f:$ $[0,1]\to[0,1]$ be an increasing function with $f(0)=0$, $f(1)=1$. Then $\mu=f(P)$ is a capacity on $\mathcal{F}$, called the distorted probability. Especially if $f$ is left (right) continuous then $\mu$ is continuous from below (above); if $f$ is strictly increasing on point $1$,
then $\mu$ satisfies $\mu(A\cap B)=1$ for $\mu(A)=1$ and $\mu(B)=1$.
\end{example}

\begin{corollary}\label{co1}
Let $(V,v)$ be a pair of continuous upper and lower probabilities  on $\mathcal{F}$ with $v(\mathcal{F}_0)=\{0,1\}$, then for any $\mathcal{F}_0$-measurable random variable $\xi$, we have the following equalities:
\begin{equation}\label{eq4}
v\left(\left\{\omega:\ \int_{\Omega}\xi dv \le\xi(\omega)\le\int_{\Omega}\xi dV\right\}\right)=1
\end{equation}
\begin{equation}\label{eq5}
V\left(\left\{\omega:\ \xi(\omega)=\int_{\Omega}\xi dv\right\}\right)=1
\end{equation}
\begin{equation}\label{eq6}
V\left(\left\{\omega:\ \xi(\omega)=\int_{\Omega}\xi dV\right\}\right)=1.
\end{equation}
\end{corollary}
\begin{proof}
Equality (\ref{eq4}) is directly from equality (\ref{eq3}) since the lower probability $v$ satisfies $v(A\cap B)=1$ if $v(A)=1$ and $v(B)=1$.

Applying the result of Theorem \ref{le1} (i) to $v$ and $V$, we can get the following four equalities

\begin{eqnarray}
\label{eq8}
v\left(\left\{\omega:\ \xi(\omega)\ge \int_{\Omega}\xi dv \right\}\right)=1,\\
\label{eq9}
v\left(\left\{\omega:\ \xi(\omega)\le\int_{\Omega}\xi dV\right\}\right)=1,\\
\label{eq10}
V\left(\left\{\omega:\  \xi(\omega)\ge\int_{\Omega}\xi dV\right\}\right)=1,\\
\label{eq11}
V\left(\left\{\omega:\ \xi(\omega)\le\int_{\Omega}\xi dv\right\}\right)=1.
\end{eqnarray}
It is easy to see that $V(A\cap B)=1$ if $V(A)=1$ and $v(B)=1$. So (\ref{eq5}) can be deduced from (\ref{eq8}) and (\ref{eq11}) while (\ref{eq6}) can be deduced from (\ref{eq9}) and (\ref{eq10}).
\end{proof}

\begin{theorem}\label{th1}
Let $\mu$ be a continuous capacity on $\mathcal{F}$ with $\mu(\mathcal{G})=\{0,1\}$, where $\mathcal{G}$ is the set of all invariant sets under $\theta$. For any real valued random variables $\xi$, there exist $\mathcal{G}$-measurable random variables $\xi^*$ and $\xi_*$ such that
\begin{equation}\label{eq12}
\mu\left(\left\{\omega:\  \liminf\limits_{n\to\infty}{1\over n}\sum_{k=0}^{n-1}\xi(\theta^{k}(\omega))\ge \int_{\Omega}\xi_* d\mu\right\}\right)=1,
\end{equation}
and
\begin{equation}\label{eq13}
\mu\left(\left\{\omega:\  \limsup\limits_{n\to\infty}{1\over n}\sum_{k=0}^{n-1}\xi(\theta^{k}(\omega))\le\int_{\Omega}\xi^* d\bar{\mu}\right\}\right)=1.
\end{equation}
Moreover, if $\mu$ satisfies $\mu(A\cap B)=1$ for any $A,B\in\mathcal{G}$ with $\mu(A)=1$ and $\mu(B)=1$,
then
\begin{equation}\label{eq14}
\mu\left(\left\{\omega:\  \int_{\Omega}\xi_* d\mu\le\liminf\limits_{n\to\infty}{1\over n}\sum_{k=0}^{n-1}\xi(\theta^{k}(\omega))\le \limsup\limits_{n\to\infty}{1\over n}\sum_{k=0}^{n-1}\xi(\theta^{k}(\omega))\le\int_{\Omega}\xi^* d\bar{\mu}\right\}\right)=1.
\end{equation}
\end{theorem}
\begin{proof}
Let $\xi_*=\liminf\limits_{n\to\infty}{1\over n}\sum_{k=0}^{n-1}\xi(\theta^{k}(\omega))$, $\xi^*=\limsup\limits_{n\to\infty}{1\over n}\sum_{k=0}^{n-1}\xi(\theta^{k}(\omega))$. Notice that
\begin{eqnarray*}
\xi_*(\theta(\omega))&=&\liminf\limits_{n\to\infty}{1\over n}\sum_{k=0}^{n-1}\xi(\theta^{k+1}(\omega))\\
&=&\liminf\limits_{n\to\infty}{1\over n}\left[\sum_{k=0}^{n}\xi(\theta^{k}(\omega))-\xi(\omega)\right]\\
&=&\liminf\limits_{n\to\infty}\left[{n+1\over n}{1\over n+1}\sum_{k=0}^{n}\xi(\theta^{k}(\omega))-{1\over n}\xi(\omega)\right]=\xi_*(\omega),
\end{eqnarray*}
we have $\xi_*$ is $\mathcal{G}$-measurable. Similarly, $\xi^*$ is $\mathcal{G}$-measurable. Therefore (\ref{eq12})  (\ref{eq13}) and (\ref{eq14}) can be derived directly from Theorem \ref{le1}.
\end{proof}

As mentioned in the introduction, we do not think $\mu(\mathcal{G})=\{0,1\}$ used in \cite{cerreia} is adequate to define the ergodicity in a capacity space.
We will use the following example to illustrate the reason.

\begin{example}\label{ex3}
Let $\Omega=\{\omega_1,\omega_2,\omega_3,\omega_4\}$, $\mathcal{F}$ be all subsets of $\Omega$.  Define $\theta:\ \Omega\to\Omega$ by $$\theta(\omega_1)=\omega_2, \quad \theta(\omega_2)=\omega_1,\quad \theta(\omega_3)=\omega_4,\quad\theta(\omega_4)=\omega_3.$$

Let $\mu_1$ be a capacity on  $\mathcal{F}$, with $\mu_1(A)=0$ if $|A|\le 1 $, $\mu_1(A)=1$ if $|A|\ge2$, where $|A|$ denotes the  number of elements in $A$. Then $\mu_1$ is neither subadditive nor superadditive.

Let $P_1$ and $P_2$ be probabilities on $\mathcal{F}$ with
$P_1(\omega_1)=P_1(\omega_2)={1\over 2}$, $P_1(\omega_3)=P_1(\omega_4)=0$,
$P_2(\omega_1)=P_2(\omega_2)=0$,  $P_2(\omega_3)=P_2(\omega_4)={1\over 2}$.
For any $A\in\mathcal{F}$, let $\mu_2(A)=\max_{i=1,2}P_i(A)$ be an upper probability.

Then it is easy to check that $\theta$ preserves both $\mu_1$ and $\mu_2$ and the set of all invariant sets is $\mathcal{G}=\{\Omega, \emptyset, \{\omega_1,\omega_2\}, \{\omega_3,\omega_4\}\}$. Here $\mu_i(\mathcal{G})=\{0,1\}$, $i=1,2$, so they together with their conjugate capacities satisfy the ergodicity definition in \cite{cerreia}. Note the conjugate capacity of $\mu_2$ is a lower probability.
But under these two different capacities, $\Omega$ can be split into two invariant sets $\{\omega_1,\omega_2\}$ and $ \{\omega_3,\omega_4\}$ with $\mu_i(\{\omega_1,\omega_2\})=1$  and $\mu_i(\{\omega_3,\omega_4\})=1$, $i=1,2$. That is to say $\Omega$ is decomposable under $\theta$.
\end{example}
Now we give the definition of ergodic transformation in a capacity space.
\begin{definition}\label{de2}
A measurable capacity preserving transformation $\theta$ on the capacity space $(\Omega,\mathcal{F},\mu)$ is said to be ergodic (with respect to  $\mu$) if for any $\theta$-invariant set $B$ the following two conditions hold:\\
(i) $\mu(B)=0$ or $\mu(B)=1$, \\
(ii) $\mu(B)=0$ or $\mu(B^c)=0$.
\end{definition}

\begin{remark}
If $\theta$ is not ergodic with respect to  capacity $\mu$ then the space $\Omega$ can be split into two $\theta$-invariant sets $B$ and $B^c$ either each having positive capacity or one of them having positive capacity which is less than $1$.
This is to say $\theta$ is not ``irreducible''.
\end{remark}
The following example shows why we do not only consider (ii) in Definition \ref{de2} to define ergodicity.
\begin{example}\label{ex4}
Let $\Omega=\{\omega_1,\omega_2,\omega_3\}$, $\mathcal{F}$ be all subsets of $\Omega$ and $P_1,P_2,P_3$ be probabilities on  $\mathcal{F}$ with
$P_1(\omega_1)=0$, $P_1(\omega_2)= P_1(\omega_3)={1\over 2},$
$P_2(\omega_2)=0$, $ P_2(\omega_1)= P_2(\omega_3)={1\over 2},$
$P_3(\omega_3)=0$, $ P_3(\omega_1)= P_3(\omega_2)={1\over 2}.$
Let $v(A)=\min_{i=1,2,3}P_i(A)$, for any $A\in\mathcal{F}$, be a lower probablity. We consider the following transformation $\theta$ with
$$\theta(\omega_1)=\omega_2,\quad \theta(\omega_2)=\omega_1,\quad \theta(\omega_3)=\omega_3.$$
It is easy to check that $v$ is $\theta$-invariant and the set of all invariant sets with respect to  $\theta$ is $\mathcal{G}=\{\Omega,\emptyset,\{\omega_1,\omega_2\},\{\omega_3\}\}$. So for any set $B\in\mathcal{G}$, $v(B)=0$ or $v(B^c)=0$. But $v(\mathcal{G})=\{0,{1\over 2},1\}$,
and $\Omega$ can be split into two invariant sets $\{\omega_1,\omega_2\}$ and $ \{\omega_3\}$ with $v(\{\omega_1,\omega_2\})={1\over2}$.
\end{example}


It is easy to check that the following proposition holds.
\begin{proposition}\label{pr3}
Let $(V,v)$ be a pair of upper and lower probabilities on $(\Omega,\mathcal{F})$.

(i) A transformation $\theta$ being ergodic with respect to  $V$ is equivalent to
for any $\theta$-invariant set $B$, either $V(B)=0$ or $V(B^c)=0$.

(ii)  The ergodicity of  $\theta$ with respect to the upper probability $V$ implies the ergodicity of  $\theta$ with respect to the lower probability $v$.
\end{proposition}


\begin{remark}
For a lower probability $v$, it is easy to see that $v(B)=0$ or $v(B)=1$ implies $v(B)=0$ or $v(B^c)=0$. Thus,
condition (i) in Definition \ref{de2} is adequate to guarantee the ergodicity of $v$. Thus the definition of ergodicity given in \cite{cerreia}, though is not enough for a general capacity, but agrees with our definition in the case of lower probability. However, as $v(B)=0$ does not imply $V(B)=0$, thus the ergodicity of $\theta$ under the lower probability does not imply the ergodicity of $\theta$ with respect to the upper probability.

\end{remark}

Motivated by Theorems 2.6 and 2.7 in \cite{fengzhao}, we derive the following Theorems \ref{th4} and \ref{th5} as the characterizations of
ergodicity in an upper probability space.
\begin{theorem}\label{th4}
Let $V$ be an upper probability  on $(\Omega,\mathcal{F})$ with continuity from below and $\theta$ be a measurable transformation from $\Omega$ to $\Omega$ preserving $V$.
Then the following four statements: \\
(i) the transformation $\theta$ is ergodic;\\
(ii) if every $B\in\mathcal{F}$ with $V(\theta^{-1}B\triangle B)=0$, then $V(B)=0$ or $V(B^c)=0$;\\
(iii) for every $A\in \mathcal{F}$ with $V(A)>0$, we have $V\left(\left(\bigcup\limits_{n=1}^{\infty}\theta^{-n}A\right)^{c}\right)=0$;\\
(iv) for every $A,B\in \mathcal{F}$ with $V(A)>0$ and $V(B)>0$, there exists $n\in\mathbb{N}$ such that $V(\theta^{-n}A\cap B)>0$,\\
have the following relations:  (i) and (ii) are equivalent; (iii) implies (iv); (iv) implies (i).
Moreover, if $V$ is continuous, then (ii) implies (iii) and all the above four statements are
equivalent.
\end{theorem}

\begin{proof}
(ii) $\Rightarrow$ (i) is obvious from Proposition \ref{pr3}.

(i) $\Rightarrow$ (ii). Let $B\in\mathcal{F}$ with $V(\theta^{-1}B\triangle B)=0$.
Since for any $n\in\mathbb{N}$
\begin{eqnarray*}
\theta^{-n} B\triangle B&\subseteq&\bigcup_{k=0}^{n-1}(\theta^{-(k+1)} B\triangle \theta^{-k} B)
=\bigcup_{k=0}^{n-1}\theta^{-k}(\theta^{-1} B\triangle  B)
\end{eqnarray*}
then by the monotonicity and subadditivity and $\theta$-invariance of $V$,
$$
V\left(\theta^{-n} B\triangle B\right)\le V\left( \bigcup_{k=0}^{n-1}\theta^{-k}(\theta^{-1} B\triangle  B)\right)
\le \sum_{k=0}^{n-1}V\left(\theta^{-k}(\theta^{-1} B\triangle  B)\right)
= \sum_{k=0}^{n-1}V\left(\theta^{-1} B\triangle  B\right).
$$
Because of $V(\theta^{-1}B\triangle B)=0$, we have
\begin{equation}\label{eq20}
V\left(\theta^{-n} B\triangle B\right)=0.
\end{equation}
Moreover,
$$\left(\bigcup_{k=0}^{\infty}\theta^{-k}B\right)\triangle B\subseteq \bigcup_{k=0}^{\infty}\left(\theta^{-k}B\triangle B\right),$$
thus from the monotonicity of $V$ and (\ref{eq20}) we have
$$
V\left(\left(\bigcup_{k=n}^{\infty}\theta^{-k}B\right)\triangle B\right)\le V\left( \bigcup_{k=0}^{\infty}\left(\theta^{-k}B\triangle B\right)\right)
\le\sum_{k=0}^{\infty}V\left( \theta^{-k}B\triangle B\right)
=0.
$$
Immediately, we have
\begin{equation}\label{eq21}
V\left(\left(\bigcup_{k=n}^{\infty}\theta^{-k}B\right)\backslash B\right)=0,
\end{equation}
and
\begin{equation}\label{eq22}
V\left(B\backslash \left(\bigcup_{k=n}^{\infty}\theta^{-k}B\right)\right)=0.
\end{equation}
Define $B_{\infty}=\bigcap\limits_{n=0}^{\infty}\bigcup\limits_{k=n}^{\infty}\theta^{-k} B$.
Combining (\ref{eq21}), it is directly from $B_{\infty}\backslash B\subseteq\left(\bigcup_{k=n}^{\infty}\theta^{-k}B\right)\backslash B$ and the monotonicity of $V$ that
\begin{equation}\label{eq23}
V\left(B_{\infty}\backslash B\right)=0.
\end{equation}
Meanwhile,
$$B\backslash \left(\bigcup_{k=n}^{\infty}\theta^{-k}B\right) \ \uparrow \ B\backslash \bigcap_{n=1}^{\infty}\left(\bigcup_{k=n}^{\infty}\theta^{-k}B\right)=B\backslash B_{\infty},$$
by the continuity from below of $V$ and (\ref{eq22}), we have
\begin{equation}\label{eq24}
V\left(B\backslash B_{\infty}\right)=0.
\end{equation}

On the other hand, $B_{\infty}$ is an invariant set since
$$\theta^{-1} B_{\infty}=\bigcap\limits_{n=0}^{\infty}\bigcup\limits_{k=n+1}^{\infty}\theta^{-k} B=B_{\infty}.$$
By the ergodicity assumption of $V$, we have
$V(B_{\infty})=0$ or $ V(B_{\infty}^c)=0.$

If $V(B_{\infty})=0$, then by the subadditivity of $V$ and (\ref{eq24}), we have $V(B)=0$ since
$$
V(B)=V(B)-V(B_{\infty})
\le V(B)-V(B\cap B_{\infty})
\le V(B\backslash (B\cap B_{\infty}))
= V(B\backslash B_{\infty})=0.
$$

If $V(B_{\infty}^c)=0$, then similarly by the subadditivity of $V$ and (\ref{eq23}), we have $V(B^c)=0$ since
\begin{eqnarray*}
V(B^c)=V(B^c)-V(B_{\infty}^c)
&\le& V(B^c)-V(B^c\cap B_{\infty}^c)\\
&\le& V(B^c\backslash (B^c\cap B_{\infty}^c))
= V(B^c\backslash B_{\infty}^c)
= V(B_{\infty}\backslash B)=0.
\end{eqnarray*}
Hence, the statement (ii) is proved.

(iii) $\Rightarrow$ (iv). Let $A,B\in \mathcal{F}$ with $V(A)>0$ and $V(B)>0$.  From (iii), we know $V\left(\left(\bigcup\limits_{n=1}^{\infty}\theta^{-n}A\right)^{c}\right)=0$.
By the subadditivity, monotonicity and continuity from below of $V$, we have
\begin{eqnarray*}
0<V(B)
&\le& V\left(B\cap \left(\bigcup\limits_{n=1}^{\infty}\theta^{-n}A\right)\right)+V\left(B\cap \left(\bigcup\limits_{n=1}^{\infty}\theta^{-n}A\right)^c\right)\\
&\le& V\left( \bigcup\limits_{n=1}^{\infty}(B\cap\theta^{-n}A)\right)+V\left(\left(\bigcup\limits_{n=1}^{\infty}\theta^{-n}A\right)^c\right)\\
&=& V\left( \bigcup\limits_{n=1}^{\infty}(B\cap\theta^{-n}A)\right)\\
&\le& \sum\limits_{n=1}^{\infty} V\left( (B\cap\theta^{-n}A)\right).
\end{eqnarray*}
Thus there exists $n\in\mathbb{N}$ such that $V\left( (B\cap\theta^{-n}A)\right)>0$. Therefore the assertion (iv) is proved.

(iv) $\Rightarrow$ (i). Let $B$ be any invariant set. If $V(B)>0$ and $V(B^c)>0$, then by (iv) and invariant assumption of $B$, there exists $n\in\mathbb{N}$ such that
$$0<V(B^c\cap\theta^{-n}B)=V(B^c\cap B)=0$$
which derives a contradiction. Hence $V(B)=0$ or $V(B^c)=0$. Therefore by Proposition \ref{pr3}, (i) is proved.

(ii) $\Rightarrow$ (iii) under the continuity assumption of $V$. Let $A\in\mathcal{F}$ with $V(A)>0$. Define
$$A_1=\bigcup_{k=1}^{\infty}\theta^{-k}A\quad \hbox{and} \quad A_{\infty}=\bigcap_{n=1}^{\infty}\bigcup_{k=n}^{\infty}\theta^{-k}A.$$
It is easy to see that
$$\theta^{-n}A_1=\bigcup_{k=n+1}^{\infty}\theta^{-k}A\ \downarrow\  A_{\infty}.$$
It follows from the continuity  and $\theta$-invariance of $V$ that
$$V(A_{\infty})=\lim_{n\to\infty}V(\theta^{-n}A_1)=V(A_1)\ge V(\theta^{-1}A)=V(A)>0.$$
Since $\theta^{-1}A_{\infty}=A_{\infty}$, by (ii), we have
$V(A_{\infty}^c)=0$. Notice that $A_{\infty}\subseteq \bigcup_{n=1}^{\infty}\theta^{-n}A$, therefore $V((\bigcup_{n=1}^{\infty}\theta^{-n}A)^c)=0$ and (iii) is proved.

It is then obvious that all the four statements are equivalent under the continuity assumption of $V$.
\end{proof}

\begin{remark}
If $\mathcal{P}$ generating the upper probability $V$ is a subset of $\Delta^{\sigma}(\Omega,\mathcal{F})$, then $V$ is continuous from below (see Lemma 2.1 in \cite{chenwuli}).
\end{remark}

\begin{definition} (Definition 3 in \cite{denishupeng}).
In an upper probability space, we call that a statement holds quasi-surely if it holds outside a  set $A$ with $V(A)=0$.
\end{definition}

\begin{remark}
In Theorem \ref{th4}, from the equivalence between (i) and (ii), we can define the invariant sets in a wider sense on upper probability space as almost invariant set in classical ergodic theory (see \cite{Durrett}).
A set $B$ is said to be quasi invariant with respect to  $\theta$ in the upper probability space $(\Omega,\mathcal{F},V)$ if $V(\theta^{-1}B\triangle B)=0.$
If the set $B$ is quasi invariant, then $V(\theta^{-1}B\cap B)=V(B)=V(\theta^{-1}B).$
Thus we can consider the quasi invariant sets when we study the ergodicity.

The statement (iii) means that if $V(A)>0$ then $\theta^{k}\omega$ will lie into $A$ in finite steps quasi-surely.

The statement (iv) means that if $V(A)>0$ and $V(B)>0$, then those $\omega$ starting from $B$, with the flow $\theta^{k}\omega$ arriving into $A$ in finite steps have positive upper probability.
\end{remark}

\begin{theorem}\label{th5}
Let $V$ be an upper probability on $(\Omega,\mathcal{F})$ and $\theta$ be a measurable transformation from $\Omega$ to $\Omega$ preserving $V$.
Then the following three statements:\\
(i) the transformation $\theta$ is ergodic;\\
(ii) if $\xi:\Omega\to\mathbb{R}$ is bounded measurable and $\xi(\theta\cdot)=\xi(\cdot)$, then $\xi$ is a constant quasi-surely;\\
(iii) if $\xi:\Omega\to\mathbb{R}$ is measurable and $\xi(\theta\cdot)=\xi(\cdot)$ quasi-surely, then $\xi$ is a constant quasi-surely,\\
have the following relations: (iii) implies (ii); (ii) implies (i). Moreover, if $V$ is continuous from below, then (ii) is equivalent to (i). If further $V$ is continuous, then all three statements are equivalent.
\end{theorem}
\begin{proof}
The proof of (iii) $\Rightarrow$ (ii) is trivial.

We now prove (ii) $\Rightarrow$ (i).
For any invariant set $A$, $I_A(\theta\omega)=I_A(\omega)$. Thus $I_A$ is a constant quasi-surely. So $I_A = 0$ or $1$ quasi-surely.
If $I_A = 0$ quasi-surely, then $V(A) = 0$. If $I_A = 1$ quasi-surely, then $V(A^c) = 0$.  Thus $\theta$ is ergodic.

Next we prove (i) $\Rightarrow$ (iii) under the assumption that $V$ is continuous.

For any $t \in\mathbb{R}$, let $A_t = \{\omega :\  \xi(\omega)>t \}$ and $A_t^c = \{\omega :\  \xi(\omega)\le t \}$. Notice that $\theta^{-1}A_t\triangle A_t\subseteq\{\omega:\ \xi(\theta\omega)\neq\xi(\omega)\} $, we have
$V(\theta^{-1}A_t\triangle A_t)=0$ since $\xi(\theta\cdot)=\xi(\cdot)$ quasi-surely.
Since $\theta$ is ergodic and $V$ is continuous from below, by
Theorem \ref{th4}, we know that $V(A_t)=0$ or $V(A_t^c)=0$. Thus, $V(A_t)=0$ or $1$.

Let $I=\{t:\ V(A_t)=0\}$. By the continuity from above of $V$, we have
$$0=V(\{\omega:\ \xi(\omega)=\infty\})=V\left(\bigcap_{n=1}^{\infty}A_n\right)=\lim_{n\to\infty}V(A_n).$$
Thus there exists $n\in \mathbb{N}$ such that $V(A_n)=0$, that is $n\in I\neq\emptyset$. So set $t_*=\inf I$, and immediately $t_*\in I$ since $V$ is continuous from below.
Hence, for any $t\ge t_*$, we have $V( A_t)=0$ and for any $t< t_*$, we have $V( A_t)=1$
and $V( A_t^c)=0$.
Due to the continuity from below of $V$, we have $V(\{\omega:\ \xi(\omega)<t_*\})=0$.
Combining $V(\{\omega:\ \xi(\omega)>t_*\})=0$ and the subadditivity of $V$, we get $V(\{\omega:\ \xi(\omega)\neq t_*\})=0$.
Thus $\xi$ is a constant $t_*$ quasi-surely.

From the proof of (i) $\Rightarrow$ (iii) we can see the continuity from above of $V$ is used to prove $I\neq\emptyset$. But if $\xi$ is bounded, of course $I\neq\emptyset$. So we do not need the  continuity from above of $V$ assumption.
Therefore, if $V$ is continuous from below, we can prove (i) $\Rightarrow$ (ii) in the same way as (i) $\Rightarrow$ (iii).
\end{proof}

\begin{remark}
The result can be also presented in the language of transformation operator $U$ defined by $U(\xi)(\omega)=\xi(\theta\omega)$ on the space of measurable functions on $\Omega$. Theorem \ref{th5} says that an upper probability preserving map is ergodic if and only if the transformation operator $U$ has eigenvalue $1$ which is simple. In classical probability space case, see Da Prato and Zabczyk \cite{Da} or Walters \cite{Walters}.
\end{remark}

Now we give a characterization of
ergodicity  through strong law of large numbers which is a type of Birkhoff's ergodic theorem (see \cite{Birkhoff}) in an upper probability space.

\begin{theorem}\label{th7}
Let $(V,v)$ be a pair of continuous upper and lower probabilities on $\mathcal{F}$ and $\theta$ be a measurable transformation from $\Omega$ to $\Omega$ preserving $V$.
Then $\theta$ is ergodic with respect to $V$ if and only if for any bounded $\mathcal{F}$-measurable random variable $\xi$,
$$\lim\limits_{n\to\infty}\frac{1}{n}\sum_{k=0}^{n-1}\xi(\theta^{k}\cdot)$$ is a constant quasi-surely.
\end{theorem}
\begin{proof}
For any bounded $\mathcal{F}$-measurable random variable $\xi$,
let
\begin{equation}\label{eq25}
A=\left\{\omega:\ \lim_{n\to\infty}{1\over n}\sum_{k=0}^{n-1}\xi(\theta^{k}\omega) \hbox{ exists }\right\}
\end{equation}
and
\begin{equation*}
\tilde{\xi}(\omega)=\begin{cases}
\lim\limits_{n\to\infty}{1\over n}\sum\limits_{k=0}^{n-1}\xi(\theta^{k}\omega),& \text{ if }  \omega\in A\\
 0,&  \text{ if }  \omega\notin A
\end{cases}.
\end{equation*}
By Lemma \ref{pr7}, $v\left(A\right)=1.$ Therefore $\tilde{\xi}$ is a bounded $\mathcal{F}$-measurable random variable with $\tilde{\xi}(\theta\cdot)=\tilde{\xi}(\cdot)$ quasi-surely.

If $\theta$ is ergodic with respect to $V$, then by Theorem \ref{th5}, $\tilde{\xi}$ is a constant quasi-surely. Therefore $\lim\limits_{n\to\infty}\frac{1}{n}\sum_{k=0}^{n-1}\xi(\theta^{k}\cdot)$ is a constant quasi-surely.

If for any bounded $\mathcal{F}$-measurable random variable $\xi$,
$\lim\limits_{n\to\infty}\frac{1}{n}\sum_{k=0}^{n-1}\xi(\theta^{k}\cdot)$ is a constant quasi-surely.
We consider $\xi$ with $\xi(\theta\cdot)=\xi(\cdot)$, then
$\lim\limits_{n\to\infty}\frac{1}{n}\sum_{k=0}^{n-1}\xi(\theta^{k}\omega)=\xi(\omega)$. Hence $\xi$ is a constant quasi-surely. It follows from Theorem \ref{th5} that $\theta$ is ergodic with respect to $V$.
\end{proof}

Next we want to give a estimate of the ergodic average constant in Theorem \ref{th7}.

\begin{theorem}\label{th8}
Let $V$ be a  continuous and concave capacity on $\mathcal{F}$, and $v$ be the conjugate capacity to $V$. Let $\theta$ be ergodic with respect to $V$ and $\xi$ be a bounded $\mathcal{F}$-measurable random variable.
Then there exists a constant $c$ such that
$$\lim\limits_{n\to\infty}\frac{1}{n}\sum_{k=0}^{n-1}\xi(\theta^{k}\cdot)=c\hbox{ quasi-surely and }c\in\left[\int_{\Omega}\xi dv,\int_{\Omega}\xi dV\right].$$
\end{theorem}
\begin{remark}
Under the assumption of Theorem \ref{th8}, capacities $(V, v)$ actually is a pair of continuous upper and lower probabilities which can be seen from Proposition \ref{pr4}.
\end{remark}
The following is the proof of Theorem \ref{th8}.

\begin{proof}
By Theorem \ref{th7}, there exists a constant $c$ such that
$$v\left(\left\{\omega:\ \lim_{n\to\infty}{1\over n}\sum_{k=0}^{n-1}\xi(\theta^{k}(\omega)) =c\right\}\right)=1.$$
We only need to prove  $c\in[\int_{\Omega}\xi dv,\int_{\Omega}\xi dV]$.
Since $\theta$ is ergodic with respect to $V$, we have $v(\mathcal{G})=\{0,1\}$.
It follows from  Theorem \ref{th1} that
\begin{eqnarray*}
v\left(\left\{\omega:\ \limsup_{n\to\infty}\frac{1}{n}\sum_{k=0}^{n-1}\xi(\theta^{k}\omega)
\le\int_{\Omega}\limsup_{n\to\infty}\frac{1}{n}\sum_{k=0}^{n-1}\xi(\theta^{k}\omega) dV\right\}\right)=1.
\end{eqnarray*}
Let $A$ be the set defined by (\ref{eq25}) in the proof of Theorem \ref{th7}, then $A$ is a $\theta$-invariant set and $V(A^c)=0$.
Combining above two equalities, by the subadditivity of the Choquet integral with respect to $V$ (Lemma \ref{le5}), we get
\begin{eqnarray*}
c
&\le&\int_{\Omega}\left(\lim_{n\to\infty}\frac{1}{n}\sum_{k=0}^{n-1}\xi(\theta^{k}\omega)\right)I_A(\omega)dV + \int_{\Omega} \left(\limsup_{n\to\infty}\frac{1}{n}\sum_{k=0}^{n-1}\xi(\theta^{k}\omega)\right)I_{A^c}(\omega)dV\\
&=&\int_{\Omega}\left(\lim_{n\to\infty}\frac{1}{n}\sum_{k=0}^{n-1}\xi(\theta^{k}\omega)\right)I_A(\omega)dV \\
&=&\lim_{n\to\infty}\int_{\Omega}\left(\frac{1}{n}\sum_{k=0}^{n-1}\xi(\theta^{k}\omega)\right)I_A(\omega)dV \\
&\le&\lim_{n\to\infty}\frac{1}{n}\sum_{k=0}^{n-1}\int_{\Omega}\xi(\theta^{k}\omega)I_A(\theta^{k}\omega)dV \\
&=&\int_{\Omega}\xi(\omega)I_A(\omega)dV\\
&=&
\int_{\Omega}\xi(\omega)dV,
\end{eqnarray*}
where the second equality is due to the dominated convergence theorem (Lemma \ref{le2}) and
the  penultimate equality is come from $\theta$ preserving the Choquet integral with respect to $V$ (Proposition \ref{pr2}).
Now we consider $-\xi$ and  $-c\le \int_{\Omega}-\xi(\omega)dV$, that is $c\ge \int_{\Omega}\xi(\omega)dv$. The proof of this theorem is completed.
\end{proof}

\begin{remark}\label{re2}
Cerreia-Vioglio, Maccheroni and Marinacci in \cite{cerreia} obtained that for a  $\theta$-invariant continuous lower probability $v$ and any bounded $\mathcal{F}$-measurable random variable $\xi$,  $\lim\limits_{n\to\infty}{1\over n}\sum_{k=0}^{n-1}\xi(\theta^{k}\cdot)$ exists
quasi-surely. If $\theta$ is ergodic in their sense,  they only showed that the random variable $\lim\limits_{n\to\infty}{1\over n}\sum_{k=0}^{n-1}\xi(\theta^{k}\cdot)$ lies in the interval
$$\left[\int_{\Omega} \limsup\limits_{n\to\infty}{1\over n}\sum_{k=0}^{n-1}\xi(\theta^{k}(\omega)) dv,\int_{\Omega}\limsup\limits_{n\to\infty}{1\over n}\sum_{k=0}^{n-1}\xi(\theta^{k}(\omega))dV\right]$$
quasi-surely. And then, for a continuous convex  capacity $v$, if it is further strongly invariant that requires all the probabilities in core$(v)$ need to be $\theta$-invariant, they confirmed the interval is $[\int_{\Omega}\xi dv,\int_{\Omega}\xi dV]$. Our result in Theorem \ref{th7} says that the ergodicity of $\theta$ with respect to $V$ is equivalent to $\lim\limits_{n\to\infty}{1\over n}\sum_{k=0}^{n-1}\xi(\theta^{k}\cdot)$ being a constant quasi-surely. This means Birkhoff's strong law of large numbers is a necessary and sufficient condition for ergodicity.
Moreover, for a convex continuous capacity $v$, we prove the constant lies in $[\int_{\Omega}\xi dv,\int_{\Omega}\xi dV]$. For this result we do not need the assumption of $v$ being strongly invariant.
\end{remark}

\section{Ergodicity of stationary processes on capacity spaces}

The notion of  stationary stochastic process on a capacity space  generalizing the usual notion of stationary stochastic process in classical probability space, was given in \cite{cerreia} as follows.

\begin{definition}\label{de3}
Given a capacity space $(\Omega,\mathcal{F},\mu)$, we say that stochastic process $\{Y_n\}_{n\in\mathbb{N}}$ is stationary if and
only if for each $n\in\mathbb{N}$, $k\in\mathbb{N}_0$ and Borel subset $A$ of $\mathbb{R}^{k+1}$,
$$
\mu(\{\omega:\ (Y_n(\omega),\cdots,Y_{n+k}(\omega))\in A\})
=\mu(\{\omega:\ (Y_{n+1}(\omega),\cdots,Y_{n+1+k}(\omega))\in A\}).
$$
\end{definition}
It is easy to see that $\{Y_n\}_{n\in\mathbb{N}}$ is stationary on $(\Omega,\mathcal{F},\mu)$ if and only if $\{Y_n\}_{n\in\mathbb{N}}$ is stationary on the conjugate capacity space $(\Omega,\mathcal{F},\overline{\mu})$
In classical probability theory, independent identically distributed random variables sequence must be stationary. However, such a result will not be valid in the capacity theory, we can find such a case in Example \ref{ex2}.

Let $(\mathbb{R}^{\mathbb{N}},\sigma(\mathcal{C}))$ denote the space of sequences endowed with the $\sigma$-algebra generated
by the set of all cylinders $\mathcal{C}$. We denote a generic element of $\mathbb{R}^{\mathbb{N}}$ by $\textbf{x}$. Any set $C$ in $\mathcal{C}$ called cylinder, has the following form
\begin{equation}\label{eq18}
C=\{\textbf{x}=(x_1,x_2,x_3,\cdots):\ (x_1,\cdots,x_n)\in H\}
\end{equation}
where $n\in\mathbb{N}$ and $H\in\mathcal{B}(\mathbb{R}^{n})$.
It is well known that $\mathcal{C} $ is an algebra.
We consider the shift transformation $\tau : \mathbb{R}^{\mathbb{N}}\to \mathbb{R}^{\mathbb{N}}$ defined by
$$\tau(\textbf{x})=\tau(x_1,x_2,x_3,\cdots)=(x_2,x_3,x_4,\cdots),\quad \hbox{for any } \textbf{x}=(x_1,x_2,x_3,\cdots)\in \mathbb{R}^{\mathbb{N}}.$$
The stochastic process $\{Y_n\}_{n\in\mathbb{N}}$ induces a measurable map from $(\Omega,\mathcal{F})$ to $(\mathbb{R}^{\mathbb{N}},\sigma(\mathcal{C}))$ by
$$\omega\mapsto\ \textbf{Y}(\omega)=(Y_1(\omega),Y_2(\omega),Y_3(\omega),\cdots),\quad \hbox{for any }\omega\in\Omega.$$
Define $\mu_{\textbf{Y}}:\sigma(\mathcal{C})\to[0,1]$ by
$$\mu_{\textbf{Y}}(C)=\mu(\textbf{Y}^{-1}(C)),\quad\hbox{for any } C\in\sigma\mathcal(C).$$
It is easy to check that $\mu_{\textbf{Y}}$ is a capacity on $\sigma(\mathcal{C})$ and $\mu_{\textbf{Y}}$ is continuous/convex/ concave if $\mu$ is continuous/convex/concave respectively,  as $\textbf{Y}^{-1}\left(\bigcup\limits_{n=1}^{\infty}C_n\right)=\bigcup\limits_{n=1}^{\infty}\textbf{Y}^{-1}\left(C_n\right)$ and $\textbf{Y}^{-1}\left(\bigcap\limits_{n=1}^{\infty}C_n\right)=\bigcap\limits_{n=1}^{\infty}\textbf{Y}^{-1}\left(C_n\right)$, for any $\{C_n\}_{n\in\mathbb{N}}\subseteq\sigma(\mathcal{C})$. Moreover, $\overline{\mu_{\textbf{Y}}}=\bar{\mu}_{\textbf{Y}}$.

\begin{proposition}\label{le3}
Let $\textbf{Y}=\{Y_n\}_{n\in\mathbb{N}}$ be a stochastic process on the capacity space $(\Omega,\mathcal{F},\mu)$ and $\mu$ is continuous.
Then $\textbf{Y}=\{Y_n\}_{n\in\mathbb{N}}$ is stationary if and only if $\mu_{\textbf{Y}}$ is the shift transformation $\tau$-invariant.
\end{proposition}
\begin{proof}
On one hand, assume $\textbf{Y}=\{Y_n\}_{n\in\mathbb{N}}$ is stationary
and let $$\mathcal{M}=\{A:\ A\in\sigma(\mathcal{C}), \mu_{\textbf{Y}}(\tau^{-1}(A))=\mu_{\textbf{Y}}(A)\}.$$
Then by the stationarity of $\textbf{Y}$, for any $C\in \mathcal{C}$ with $H$ given in (\ref{eq18}) corresponding to set $C$, we have
\begin{eqnarray*}
\mu_{\textbf{Y}}(\tau^{-1}(C))=\mu(\{\omega:\ \tau \textbf{Y}(\omega)\in C\})
&=&\mu(\{\omega:\ (Y_2(\omega),\cdots,Y_{n+1}(\omega))\in H\})\\
&=&\mu(\{\omega:\ (Y_1(\omega),\cdots,Y_{n}(\omega))\in H\})\\
&=&\mu(\{\omega:\ \textbf{Y}(\omega)\in C\})=\mu_{\textbf{Y}}(C),
\end{eqnarray*}
which implies that $\mathcal{C}\subseteq\mathcal{M}$.
Because of the continuity of $\mu$, $\mu_{\textbf{Y}}$ is continuous and then $\mathcal{M}$ is a monotone class.
Therefore, by the monotone class theorem we can get $\mathcal{M}=\sigma(\mathcal{C})$ which means $\mu_{\textbf{Y}}$ is the shift transformation $\tau$-invariant.

On the other hand, assume that $\mu_{\textbf{Y}}$ is the shift transformation $\tau$-invariant, then
 for each $n\in\mathbb{N}$, for each $k\in\mathbb{N}_0$, and for each $A\in\mathcal{B}(\mathbb{R}^{k+1})$,
\begin{eqnarray*}
\mu(\{\omega:\ (Y_n(\omega),\cdots,Y_{n+k}(\omega))\in A\})
&=& \mu_{\textbf{Y}}(\mathbb{R}^{n-1}\times A\times \mathbb{R}^{\mathbb{N}-n-k})\\
&=& \mu_{\textbf{Y}}(\tau^{-1}(\mathbb{R}^{n-1}\times A\times \mathbb{R}^{\mathbb{N}-n-k}))\\
&=& \mu_{\textbf{Y}}(\mathbb{R}^{n}\times A\times \mathbb{R}^{\mathbb{N}-n-k-1})\\
&=&\mu(\{\omega:\ (Y_{n+1}(\omega),\cdots,Y_{n+1+k}(\omega))\in A\}).
\end{eqnarray*}
That is to say $\{Y_n\}_{n\in\mathbb{N}}$ is stationary.
\end{proof}

\begin{definition}\label{de3}
The stochastic process $\{Y_n\}_{n\in\mathbb{N}}$ on capacity space $(\Omega,\mathcal{F},\mu)$  is called ergodic if the shift transformation $\tau$ is ergodic with respect to $\mu_{\textbf{Y}}$.
\end{definition}

Now we give the strong law of large numbers for stationary and ergodic stochastic sequences on a capacity space.
\begin{theorem}\label{th2}
Let $(V,v)$ be a pair of continuous upper and lower probabilities on $\mathcal{F}$.
If a bounded stationary process $\textbf{Y}=\{Y_n\}_{n\in\mathbb{N}}$ on capacity space $(\Omega,\mathcal{F},V)$  is ergodic, then
 there exists a constant $c$ such that
$$
v\left(\left\{\omega\in\Omega:\ \lim\limits_{n\to\infty}{1\over n}\sum_{k=1}^{n}Y_k(\omega)=c\right\}
\right)=1.
$$
If further $V$ is concave, that is $v$ is convex, then $c\in [\int_{\Omega}Y_1d v,\int_{\Omega}Y_1d V]$.
\end{theorem}
\begin{proof}
It is easy to check that $(V_\textbf{Y},v_\textbf{Y})$ is a pair of  continuous upper and lower probabilities on $\sigma(\mathcal{C})$ since $(V,v)$ is a pair of continuous upper and lower probabilities on $\mathcal{F}$.
By Proposition \ref{le3},
$v_{\textbf{Y}}$ and $V_{\textbf{Y}}$ are the shift transformation $\tau$-invariant. Define $\xi:\mathbb{R}^{\mathbb{N}}\to\mathbb{R}$ by $\xi(\textbf{x})=\xi(x_1,x_2,x_3,\cdots)=x_1$, for any $\textbf{x}=(x_1,x_2,x_3,\cdots)\in\mathbb{R}^{\mathbb{N}}$.
Since $\tau$ is ergodic with respect to $V_{\textbf{Y}}$,
then we can get the following equality by Theorem \ref{th7} that there exists a constant $c$ such that
$$v_\textbf{Y}\left(\left\{\textbf{x}: \lim\limits_{n\to\infty}{1\over n}\sum_{k=1}^{n}\xi(\tau^{k-1}\textbf{x})=c\right\}\right)=1. $$
Notice ${1\over n}\sum_{k=1}^{n}x_k={1\over n}\sum_{k=1}^{n}\xi(\tau^{k-1}\textbf{x})$,
we have
$$
1=v_\textbf{Y}\left(\left\{\textbf{x}:\ \lim\limits_{n\to\infty}{1\over n}\sum_{k=1}^{n}x_k=c\right\}\right)
=v\left(\left\{\omega:\ \lim\limits_{n\to\infty}{1\over n}\sum_{k=1}^{n}Y_k(\omega)=c\right\}\right).
$$
If further $V$ is concave, then $V_\textbf{Y}$ is concave while $v$ and $v_\textbf{Y}$ are convex.
Thus, it follows from Theorem \ref{th8} that $c\in[\int_{\mathbb{R^{\mathbb{N}}}}\xi dv_\textbf{Y},\int_{\mathbb{R^{\mathbb{N}}}}\xi dV_\textbf{Y}]$.
By the transformation rule of Choquet integral (see Proposition 5.1 in \cite{Denneberg}), we have
$$\int_{\mathbb{R^{\mathbb{N}}}}\xi dv_\textbf{Y}=\int_{\mathbb{R^{\mathbb{N}}}}\xi dv(\textbf{Y}^{-1})=\int_{\Omega}\xi(\textbf{Y}) dv=\int_{\Omega}Y_1 dv$$
and similarly
$$\int_{\mathbb{R^{\mathbb{N}}}}\xi dV_\textbf{Y}=\int_{\Omega}Y_1 dV.$$
As a consequence,  we completed the proof of Theorem \ref{th2}.
\end{proof}

\begin{corollary}\label{co2}
Let $(V,v)$ be a pair of continuous upper and lower probabilities on $\mathcal{F}$.
If a bounded stochastic process $\textbf{Y}=\{Y_n\}_{n\in\mathbb{N}}$ on capacity space $(\Omega,\mathcal{F},V)$  is stationary and for any $n\in\mathbb{N}$, $\sigma(Y_k,\ k\le n)$
and  $\sigma(Y_k,\ k\ge n+1)$ are independent with respect to  $V$.
Then there exists a constant $c$ such that
$$
v\left(\left\{\omega\in\Omega:\ \lim\limits_{n\to\infty}{1\over n}\sum_{k=1}^{n}Y_k(\omega)=c\right\}
\right)=1.
$$
If further $V$ is concave, that is $v$ is convex, then $c\in [\int_{\Omega}Y_1d v,\int_{\Omega}Y_1d V]$.
\end{corollary}
\begin{proof}
From Theorem \ref{th2}, we only need to prove that the shift transformation $\tau$ is ergodic  with respect to $V_{\textbf{Y}}$.
Let $B$ be any $\tau$-invariant set, then for any $n\in\mathbb{N}$
\begin{eqnarray*}
\textbf{Y}^{-1}(B)&=&\textbf{Y}^{-1}(\tau^{-1}(B))\\
&=&\textbf{Y}^{-1}(\tau^{-n}(B))\\
&=&\{\omega:\ \tau^{n}\textbf{Y}(\omega)\in B\}\\
&=&\{\omega:\ (Y_{n+1}(\omega),Y_{n+2}(\omega),\cdots)\in B\}\in\sigma(Y_k,\ k\ge n+1).
\end{eqnarray*}
So $\textbf{Y}^{-1}(B)\in\mathcal{T}$. By the assumption that  for any $n\in\mathbb{N}$, $\sigma(Y_k,\ k\le n)$
and  $\sigma(Y_k,\ k\ge n+1)$ are independent with respect to  $V$, it is follows from the Kolmogorov $0$-$1$ Law in capacity spaces (Theorem \ref{th3}) that $V(\textbf{Y}^{-1}(B))=0$ or $V(\textbf{Y}^{-1}(B)^c)=0$.
Therefore $V_{\textbf{Y}}(B)=0$ or $V_{\textbf{Y}}(B^c)=0$, that is $\tau$ is ergodic  with respect to $V_{\textbf{Y}}$ by Proposition \ref{pr3}.
\end{proof}

\begin{remark}

There are many references on strong law of large numbers for capacities under different definitions of independence and identical distributions, see for example   \cite{cerreia}, \cite{chen}, \cite{chenwuli}, \cite{Cooman1}, \cite{epstein}, \cite{Marinacci1}, \cite{Marinacci2}, \cite{peng6}, \cite{Teran}, \cite{zhang}  and references therein.
Comparing with these papers, we replace the independent identically distributed hypothesis  by the stationarity and ergodicity. We weaken the assumption of total monotonicity of $v$  in \cite{Marinacci1} and \cite{Marinacci2} to convexity while we need the continuity of the lower probability $v$.
But we do not need $\Omega$ to be a Polish space as in \cite{Marinacci1},  or compact space as in \cite{Marinacci2} or  a finite space as in \cite{epstein}. It
was obtained in \cite{cerreia} that the empirical average exists quasi-surely. By improving the definition of ergodicity,
we obtain that the empirical average is constant quasi-surely. This is a property that was not present in previous work.

\end{remark}

\end{document}